\documentclass[10pt,reqno]{amsart}

\usepackage{amsmath,amssymb,amsthm,fancyhdr,mathrsfs,times,bbm}
\usepackage[latin1]{inputenc}

\makeindex

\newcommand{\diff}{\mt{d}}

\newcommand{\mt}[1]{{\text{\rm #1}}}  

\newcommand{\comment}[1]{}            

\newcommand{\set}[1]{\{#1\}}
\newcommand{\bigset}[1]{\big\{#1\big\}}

\newcommand{\suchthat}{\;|\;}         
\newcommand{\bigsuchthat}{\;\big|\;}

\newcommand{\restrict}{|}           
\newcommand{\compose}{\circ}          
\newcommand{\define}{\;{\rm :=}\;}

\newcommand{\without}{\mathord{\setminus}}

\newcommand{\id}{\mt{id}}             
\newcommand{\pr}{\mt{pr}}             
\newcommand{\R}{\mathbb{R}}           
\newcommand{\N}{\mathbb{N}}           
\newcommand{\Z}{\mathbb{Z}}           
\newcommand{\leer}{\varnothing}       
\newcommand{\mfbd}{\partial}          
\newcommand{\eval}[2]{\langle #1,#2 \rangle}     
\newcommand{\bigeval}[2]{\big\langle #1,#2 \big\rangle}

\newcommand{\abs}[1]{\lvert#1\rvert}  
\newcommand{\bigabs}[1]{\big\lvert#1\big\rvert}

\newcommand{\eps}{\varepsilon}        
\newcommand{\oointerval}[2]{\left]#1,#2\right[}  
\newcommand{\cointerval}[2]{\left[#1,#2\right[}  

\makeatletter
\def\moverlay{\mathpalette\mov@rlay}
\def\mov@rlay#1#2{\leavevmode\vtop{%
\baselineskip\z@skip \lineskiplimit-\maxdimen
\ialign{\hfil$#1##$\hfil\cr#2\crcr}}}
\makeatother

\swapnumbers
\theoremstyle{definition}
\newtheorem{remark}{Remark}[section]
\newtheorem{fact}[remark]{Fact}
\newtheorem{definition}[remark]{Definition}
\newtheorem{notation}[remark]{Notation}
\theoremstyle{plain}
\newtheorem{theorem}[remark]{Theorem}
\newtheorem{lemma}[remark]{Lemma}

\newtheorem{proposition}[remark]{Proposition}

\newcommand{\Sym}{\mt{Sym}}       

\newcommand{\Riem}{\mathord{\mt{Riem}}}
\newcommand{\SecondFF}{\mathord{\mathit{II}}}
\newcommand{\dist}{\mathord{\mt{dist}}}
\newcommand{\length}{\mathord{\mt{length}}}

\newcommand{\eucl}{\mt{eucl}}

\renewcommand{\bar}{\overline}

\newcommand{\LeviCivita}{Levi-Civita}

\newcommand{\M}{\mathcal{M}}

\newcommand{\interior}{\textsl{interior}}
\newcommand{\exterior}{\textsl{exterior}}

\begin{document}

\title{A remark on the rigidity case of the positive energy theorem}
\author{Marc Nardmann}
\address{Department of Mathematics, University of Regensburg}
\email{Marc.Nardmann\@@mathematik.uni-regensburg.de}
\thanks{Supported by the Deutsche Forschungsgemeinschaft within the
priority programme ``Globale Differentialgeometrie''.}

\begin{abstract}
In their proof of the positive energy theorem, Schoen and Yau showed that every asymptotically flat spacelike hypersurface $M$ of a Lorentzian manifold which is flat along $M$ can be isometrically imbedded with its given second fundamental form into Minkowski spacetime as the graph of a function $\R^n\to\R$; in particular, $M$ is diffeomorphic to $\R^n$. In this short note, we give an alternative proof of this fact. The argument generalises to the asymptotically hyperbolic case, works in every dimension $n$, and does not need a spin structure.
\end{abstract}

\maketitle


\section{Introduction}

The \emph{rigidity case} of the positive energy theorem is the situation when $E=\abs{P}$ holds for the energy $E\in\R$ and the momentum $P\in\R^n$ of an asymptotically flat spacelike hypersurface $M$ of a Lorentzian $(n+1)$-manifold $(\bar{M},\bar{g})$ with $n\geq3$ which satisfies the dominant energy condition at every point of $M$. The positive energy theorem says that then the Riemann tensor of $\bar{g}$ vanishes at every point of $M$; we call this the \emph{rigidity statement}.

\smallskip
This has been proved by Parker/Taubes \cite{ParkerTaubes} in the case when $M$ admits a spin structure --- and under the assumption that $M$ is $3$-dimensional, but the argument generalises to higher dimensions. (The original proof of Witten \cite{Witten1981} made the slightly stronger assumption that $(\bar{M},\bar{g})$ satisfies the dominant energy condition on a neighbourhood of $M$.)

\smallskip
Another proof of the positive energy theorem, in particular of the rigidity statement, had been given earlier by Schoen/Yau \cite{SchoenYau1,SchoenYau2,SchoenYau3}, without the spin assumption --- again assuming $n=3$, but the argument can be generalised to $n\leq7$. More recently, Lohkamp extended their approach to higher dimensions \cite{Lohkamp2008arXiv}; the details for arbitrary fundamental forms have not been published yet, however. Schoen has announced a proof in a similar spirit.

\smallskip
Schoen/Yau proved actually more than Parker/Taubes: they showed that in the rigidity case the Riemannian $n$-manifold $M$ with its second fundamental form induced by the imbedding in $(\bar{M},\bar{g})$ can be imbedded isometrically into Minkowski spacetime $\R^{n,1}=\R^n\times\R$ as the graph of a function $\R^n\to\R$, which implies in particular that $M$ is diffeomorphic to $\R^n$.

\smallskip
It is natural to ask whether one can decouple the proof of imbeddability into Minkowski spacetime from the proof of the rigidity statement: When we know already --- for instance from the Parker/Taubes proof --- that $\bar{g}$ is flat along $M$, can we deduce directly that $M$ with its second fundamental form admits an imbedding of the desired form and is in particular diffeomorphic to $\R^n$?

\smallskip
The aim of the present short article is to show how this can be done in a simple way, independently of the Schoen/Yau arguments, and with minimal assumptions. Locally, the desired imbeddability follows already from the fundamental theorem of hypersurface theory due to Bär/Gauduchon/Moroianu \cite[Section 7]{BaerGauduchonMoroianu} (which has a short elegant proof).

\smallskip
Since this theorem applies not only to flat metrics but to metrics of arbitrary constant sectional curvature, we can also consider the case of imbeddings into anti-de Sitter spacetime. An analogue of the Parker/Taubes proof in this situation is the work by Maerten \cite{Maerten2006}, which requires a spin assumption. He shows in this case that the hypersurface with its second fundamental form imbeds isometrically into anti-de Sitter spacetime. As Schoen/Yau, he does this via an explicit construction which is a by-product of the specific method that is used to prove the positive energy theorem.

\smallskip
The result of the present article, Theorem \ref{rigid} below, applies in a situation when it has already been proved somehow that along the hypersurface the Gauss and Codazzi equations of an ambient Lorentzian metric of constant curvature $c\leq0$ are satisfied. The conclusion is that then a suitable isometric imbedding into Minkowski or anti-de Sitter spacetime exists and is essentially unique, which implies in particular that the hypersurface is diffeomorphic to $\R^n$. The proof does not require any spin assumption or dimensional restriction.

\smallskip
Let us adopt the following conventions and terminology. All manifolds, bundles, metrics, maps, etc.\ are smooth. The sign convention for the Riemann tensor is $\Riem(u,v)w = \nabla_u\nabla_vw -\nabla_v\nabla_uw -\nabla_{[u,v]}w$. Lorentzian metrics on $(n+1)$-manifolds have signature $(n,1)$ (i.e.\ $n$ positive=spacelike dimensions, $1$ negative=timelike dimension).

\begin{definition}[hypersurface data set]
A \emph{hypersurface data set} is a quadruple $(M,g,N,K)$ such that $M$ is a manifold, $g$ is a Riemannian metric on $M$, \;$N$ is a Riemannian line bundle over $M$ (i.e.\ a real line bundle equipped smoothly with scalar products on the fibres), and $K$ is a section in $\Sym^2(T^\ast M)\otimes N\to M$.

\smallskip
When $M$ is a spacelike hypersurface of a Lorentzian manifold $(\bar{M},\bar{g})$, then the \emph{hypersurface data set induced by the inclusion $M\to(\bar{M},\bar{g})$} is the hypersurface data set $(M,g,N,K)$ such that $g$ is the restriction of $\bar{g}$, such that $N$ is the normal bundle of $M$ in $(\bar{M},\bar{g})$ equipped with the restriction of $-\bar{g}$ as fibre metric, and such that $K$ is the second fundamental form of $M$ in $(\bar{M},\bar{g})$.

\smallskip
Let $(M,g,N,K)$ be a hypersurface data set. An \emph{isometric imbedding of $(M,g,N,K)$} into a Lorentzian manifold $(\bar{M},\bar{g})$ is a pair $(f,\iota)$ such that
\begin{itemize}
\item
$f\colon(M,g)\to(\bar{M},\bar{g})$ is an isometric imbedding;
\item
$\iota$ is an isomorphism of Riemannian line bundles from $N$ to the normal bundle $N'$ of the spacelike hypersurface $M'\define f(M)$ in $(\bar{M},\bar{g})$, where the fibre metric on $N'$ is the restriction of $-\bar{g}$;
\item
the second fundamental form $\SecondFF\in\Gamma(\Sym^2T^\ast M'\otimes N')$ of $M'$ in $(\bar{M},\bar{g})$ is given by $\SecondFF(f_\ast v,f_\ast w) = \iota(K(v,w))$ for all $x\in M$ and $v,w\in T_xM$.
\end{itemize}
An \emph{isometric immersion of $(M,g,N,K)$} into $(\bar{M},\bar{g})$ is a pair $(f,\iota)$ such that $f\colon M\to\bar{M}$ is an immersion, such that $\iota$ is a map whose domain is the total space of $N$, and such that every $x\in M$ has a neighbourhood $U$ for which $(f\restrict U,\iota\restrict(N\restrict U))$ is an isometric imbedding of $(U,g\restrict U,N\restrict U,K\restrict U)$ into $(\bar{M},\bar{g})$.
\end{definition}

\emph{Remark.} In most contexts where a spacelike hypersurface $M$ of a Lorentzian manifold $(\bar{M},\bar{g})$ is considered (e.g.\ in the positive energy theorem or discussions of the constraint equations in General Relativity), it is assumed that the normal bundle of $M$ is trivial (i.e.\ that $\bar{g}$ is time-orientable on a neighbourhood of $M$), and a unit normal vector field is fixed. This assumption is often unnecessary, in particular for the rigidity case of the positive energy theorem: We obtain the triviality of the normal bundle as a \emph{conclusion}, we do not have to assume it.

\begin{definition}
Let $(M,g,N,K)$ be a hypersurface data set. We denote the fibre scalar product on $N$ by $\eval{.}{.}_N$. We define a covariant derivative $\diff^N$ on the Riemannian line bundle $N\to M$ by declaring every local unit-length section to be parallel. We define $\nabla^{g,N}$ to be the covariant derivative on the vector bundle $\Sym^2T^\ast M\otimes N\to M$ induced by the {\LeviCivita} connection of $g$ and $\diff^N$.

\smallskip
Let $c\in\R$. \;$(M,g,N,K)$ \emph{satisfies the Gauss and Codazzi equations for constant curvature $c$} iff the equations
\begin{align*}
c\big(g(u,z)g(v,w) -g(u,w)g(v,z)\big) &= \Riem_g(u,v,w,z) -\eval{K(u,w)}{K(v,z)}_N +\eval{K(u,z)}{K(v,w)}_N \;\;,\\
0 &= -\bigeval{(\nabla^{g,N}_uK)(v,w)}{n}_N +\bigeval{(\nabla^{g,N}_vK)(u,w)}{n}_N
\end{align*}
hold for all $x\in M$ and $u,v,w,z\in T_xM$ and $n\in N_x$.
\end{definition}

\begin{fact}
Let $(M,g,N,K)$ be the hypersurface data set induced by the inclusion of a spacelike hypersurface $M$ into a Lorentzian manifold $(\bar{M},\bar{g})$ which has constant (sectional) curvature $c$ at every point of $M$. Then $(M,g,N,K)$ satisfies the Gauss and Codazzi equations for constant curvature $c$.\qed
\end{fact}

\emph{Remark.} When the hypersurface data set $(M,g,N,K)$ induced by the inclusion of a spacelike hypersurface $M$ into a Lorentzian manifold $(\bar{M},\bar{g})$ satisfies the Gauss and Codazzi equations for constant curvature $c$, then $(\bar{M},\bar{g})$ does in general not have constant curvature $c$ at any point of $M$. The reason is that the Gauss and Codazzi equations do not yield information about the curvature components $\Riem_{\bar{g}}(n,v,w,n)$ with $v,w\in T_xM$ and $n\in N_x$.

\begin{notation} \label{notation}
Let $n,r\geq0$, let $c\in\R_{\leq0}$. Let $\R^{n,r}$ denote $\R^{n+r}$ equipped with the semi-Riemannian metric $g_{n,r}\define \sum_{i=1}^n\diff x_i^2 -\sum_{i=n+1}^{n+r}\diff x_i^2$. We define $\M^{n,1}_0$ to be Minkowski spacetime $\R^{n,1}$. For $c<0$, we consider the pseudohyperbolic spacetime $\mathcal{H}^{n,1}_c\define \bigset{x\in\R^{n,2} \bigsuchthat g_{n,2}(x,x)=\frac{1}{c}}$ (which is a Lorentzian submanifold of $\R^{n,2}$) and its universal covering $\varpi\colon \R^n\times\R \to\mathcal{H}^{n,1}_c$ given by $(x,t)\mapsto (x,\cos t\sqrt{\abs{x}^2-1/c},\sin t\sqrt{\abs{x}^2-1/c})$, and we define the anti-de Sitter spacetime $\M^{n,1}_c$ to be $\R^n\times\R$ equipped with the $\varpi$-pullback metric of the metric on $\mathcal{H}^{n,1}_c$. (Both $\mathcal{H}^{n,1}_c$ and $\M^{n,1}_c$ have constant curvature $c$; sometimes $\mathcal{H}^{n,1}_c$ instead of $\M^{n,1}_c$ is called anti-de Sitter spacetime.)

\smallskip
For $c\leq0$, we define $\pr\colon\M^{n,1}_c=\R^n\times\R\to\R^n$ to be the projection $(x,t)\mapsto x$.
\end{notation}

Now we can state the main result (our definition of \emph{simply connected} includes being connected):

\begin{theorem} \label{rigid}
Let $n\geq0$ and $c\in\R_{\leq0}$, let $M$ be a connected $n$-manifold which contains a simply connected noncompact $n$-di\-mensional sub\-mani\-fold-with-boundary that is closed in $M$ and has compact boundary, let $(M,g,N,K)$ be a hypersurface data set which satisfies the Gauss and Codazzi equations for constant curvature $c$. Assume that $(M,g)$ is complete. Then:
\begin{enumerate}
\item
$(M,g,N,K)$ admits an isometric imbedding $(f,\iota)$ into $\M^{n,1}_c$ such that $\pr\compose f\colon M\to\R^n$ is a diffeomorphism.
\item
When $(\tilde{f},\tilde{\iota})$ is an isometric immersion of $(M,g,N,K)$ into $\M^{n,1}_c$, then there is an isometry $A\colon\M^{n,1}_c\to\M^{n,1}_c$ with $\tilde{f} = A\compose f$; in particular, $\tilde{f}$ is an imbedding.
\end{enumerate}
\end{theorem}

\emph{Remark 1.} In the rigidity case of (the asymptotically flat version of) the positive energy theorem, the assumptions of our theorem are satisfied: The hypersurface data set is induced by the inclusion of $M$ into a Lorentzian manifold which is flat along $M$, and thus satisfies the Gauss and Codazzi equations for constant curvature $0$. The Riemannian metric $g$ is complete (this follows from the definition of asymptotic flatness). $M$ contains a compact $n$-dimensional submanifold-with-boundary $C$ such that $M\without(C\without\mfbd C)$ is diffeomorphic to a nonempty disjoint union of copies of $\R^n\without\mt{(open ball)}$ each of which is closed in $M$ (this closedness follows from the completeness of the metric) and simply connected (because $n\geq3$ is assumed in the positive energy theorem).

\smallskip
Similarly, the assumptions are satisfied in Maerten's theorem for asymptotically hyperbolic hypersurfaces \cite[second half of the proof of the first theorem in Section 4]{Maerten2006}.

\medskip
\emph{Remark 2.} Statement (i) shows that $f(M)$ is the spacelike graph of a function $\R^n\to\R$. This implies also that $f(M)$ is an acausal subset of $\M^{n,1}_c$. (Note that e.g.\ not every spacelike imbedding $f\colon\R^n\to\R^{n,1}$ is acausal: consider an imbedding that winds up, i.e.\ in the direction of increasing time, in a spacelike way like a spiral staircase.)

\medskip
\emph{Remark 3.} Theorem \ref{rigid} would clearly be false without the simply-connectedness assumption, even in the case $K\equiv0$: take e.g.\ $(M,g,N,K)$ to be the hypersurface data set induced by the inclusion of $M=\R^{n-1}\times S^1\times\set{0}$ into the flat product Lorentzian manifold $\R^{n-1}\times S^1\times\R$ with $\R$ as timelike factor. Then (i) is clearly not true.

\smallskip
The theorem would also be false without the completeness assumption: small subsets (e.g.\ diffeomorphic to a ball or an annulus) of a complete spacelike hypersurface in Minkowski spacetime yield counterexamples.

\medskip
\emph{Remark 4.} The theorem does not assume that the Riemannian line bundle $N$ is trivial. But it implies that $N$ is trivial, because every Riemannian line bundle over $\R^n$ is trivial. Note that also this triviality would in general not hold without the simply-connectedness assumption: flat $\R^{n-1}\times S^1$ admits an isometric imbedding (with $K\equiv0$) into the flat Lorentzian manifold $\R^{n-1}\times\mathfrak{M}$, where $\mathfrak{M}$ is the Möbius strip, regarded as a line bundle over $S^1$ with timelike fibres. The normal bundle is not trivial in this case, but all assumptions of Theorem \ref{rigid} except for the simply-connectedness are satisfied.

\medskip
\emph{Remark 5.} $A$ in (ii) is in general neither time orientation-preserving nor space orientation-preserving. (Every isometric imbedding can be composed with an isometry of $\M^{n,1}_c$ which is space and/or time orientation-reversing.)

\medskip
\emph{Remark 6.} In the case $c<0$, the theorem holds also with $\mathcal{H}^{n,1}_c\cong \R^n\times S^1$ and the projection $\pr'\colon\R^n\times S^1\ni(x,t)\mapsto x\in\R^n$ instead of $\M^{n,1}_c$ and $\pr$. Similarly, Minkowski spacetime $\M^{n,1}_0$ is the universal cover of a Lorentzian manifold $\mathcal{H}^{n,1}_0 = (\R^n\times S^1,g_0)$ via the covering $q\colon \R^n\times\R\ni(x,s)\mapsto(x,[s])\in \R^n\times(\R/\Z)$, and the theorem would hold with $\mathcal{H}^{n,1}_0$ and $\pr'$ instead of $\M^{n,1}_0$ and $\pr$. One can see this either by checking that the proof of Theorem \ref{rigid} remains valid with these modifications, or directly by applying the theorem and composing maps $M\to\M^{n,1}_c$ with $q$.

\medskip
The rest of the article contains the proof of Theorem \ref{rigid}.


\section{The fundamental theorem for hypersurfaces}

We need the following special case of the fundamental theorem for hypersurfaces due to Bär/Gauduchon/Moroianu \cite[Section 7]{BaerGauduchonMoroianu}:

\begin{proposition} \label{bgm}
Let $n\geq0$ and $c\in\R$, let $M$ be a simply connected $n$-manifold, let $(M,g,N,K)$ be a hypersurface data set which satisfies the Gauss and Codazzi equations for constant curvature $c$. Then $(M,g,N,K)$ admits an isometric immersion into $\M^{n,1}_c$. When $f_0,f_1$ are isometric immersions of $(M,g,N,K)$ into $\M^{n,1}_c$, then there exists an isometry $A\colon\M^{n,1}_c\to\M^{n,1}_c$ with $f_1 = A\compose f_0$.
\end{proposition}
\begin{proof}[Remarks on the proof]
Bär/Gauduchon/Moroianu (BGM) consider the situation when the metric on $M$ has arbitrary signature and trivial spacelike normal bundle in $(\bar{M},\bar{g})$ (see the beginning of \cite[Section 3]{BaerGauduchonMoroianu}). Since every real line bundle over a simply connected manifold is trivial (the Stiefel/Whitney class $w_1(N)\in H^1(M;\Z_2)$ classifies real line bundles $N\to M$ up to isomorphism), so is our $N$. To apply the BGM result in our case, we reverse the signs of our $\bar{g}$ and $c$, then use their Corollary 7.5. We obtain existence, and uniqueness up to isometries, of isometric immersions of the sign-reversed version of $(M,g,N,K)$ into the sign-reversed version of $\M^{n,1}_c$. This yields existence and uniqueness up to isometries of isometric immersions of $(M,g,N,K)$ into $\M^{n,1}_c$.

\smallskip
In this argument we have not applied the BGM result literally, because the sign-reversed version of our $\M^{n,1}_c$ is the (nontrivial) universal cover of BGM's $\mathbb{M}^{1,n}_{-c}$. But the BGM Corollary 7.4, which makes only a local statement, does not care about the difference, and the BGM Corollary 7.5 then follows from a standard monodromy argument which works for every geodesically complete manifold of signature $(1,n)$ and constant curvature $-c$.
\end{proof}


\section{Quasicoverings}

Let us use the following terminology:

\begin{definition} \label{defquasicovering}
Let $M,B$ be $n$-manifolds. A map $\phi\colon M\to B$ is a \emph{quasicovering} iff it has the following properties:
\begin{enumerate}
\item
$\phi$ is an immersion (equivalently: it is a local diffeomorphism, i.e., every $y\in M$ has an open neighbourhood $U$ such that $\phi\restrict U$ is diffeomorphism onto its image).
\item
The $\phi$-preimage of every connected component of $B$ is nonempty.
\item
For all paths $\gamma\colon[0,1]\to B$ and $\tilde{\gamma}\colon\cointerval{0}{1}\to M$ with $\phi\compose\tilde{\gamma} = \gamma\restrict\cointerval{0}{1}$, there exists an extension of $\tilde{\gamma}$ to a path $[0,1]\to M$.
\end{enumerate}
\end{definition}

We will only be interested in the case $B=\R^n$.

\smallskip
It is easy to see that every covering map (in the smooth category) is a quasicovering. (Recall that a covering map is defined by the condition that every $x\in B$ has an open neighbourhood $U$ such that $\phi^{-1}(U)$ is the nonempty union of open disjoint sets $U_i$ each of which is mapped diffeomorphically onto $U$ by $\phi$.)

\smallskip
Less obviously, every quasicovering is a covering; i.e., the two concepts are equal. I do not know a reference where this elementary fact is stated explicitly, although I suspect that some exists. In the proof of Theorem \ref{rigid} below we will be in a situation where it is easy to check that a certain map $\phi\colon M\to\R^n$ is a quasicovering. If we knew a priori that it is a covering, then covering theory would imply that it is a diffeomorphism (because $\R^n$ is simply connected); this is what we need.

\smallskip
But the covering property of $\phi$ is hard to verify directly: For every $x\in B$, every $y\in\phi^{-1}(\set{x})$ has an open neighbourhood $U_y$ which is mapped diffeomorphically to an open neighbourhood $V_y$ of $x$. But $\phi^{-1}(\set{x})$ could a priori be infinite, and we would have to show that the sets $U_y$ can be chosen such that the intersection of the sets $V_y$ is a neighbourhood of $x$.

\smallskip
However, one can show directly that every quasicovering $\phi\colon M\to\R^n$ is a diffeomorphism just by going through the standard proofs of covering theory and checking that they remain valid, essentially word by word, for a quasicovering. One can even verify in this way that the classifications of coverings and quasicoverings coincide in general, which implies that every quasicovering is a covering; but we are not interested in doing that.

\begin{lemma} \label{lemmaquasicovering}
Let $M,B$ be connected $n$-manifolds with $B$ simply connected, let $\phi\colon M\to B$ be a quasicovering. Then $\phi$ is a diffeomorphism.
\end{lemma}
\begin{proof}[Sketch of proof]
As mentioned, we just have to go through some of the standard proofs of covering theory, e.g.\ as in \cite[Sections III.3--8]{Bredon}. The main steps are as follows.

\smallskip
\emph{Step 1: For every path $\gamma\colon[0,1]\to B$ and every $z\in M$ with $\phi(z)=\gamma(0)$, there exists a unique path $\tilde{\gamma}\colon[0,1]\to M$ with $\phi\compose\tilde{\gamma}=\gamma$ and $\tilde{\gamma}(0)=z$.} In order to prove this, consider the set $I$ of all $t\in[0,1]$ such that there exists a unique path $\tilde{\gamma}\colon[0,t]\to M$ with $\phi\compose\tilde{\gamma}=\gamma\restrict[0,t]$ and $\tilde{\gamma}(0)=z$. Clearly $0\in I$. Property (i) in the quasicovering definition implies that $I$ is open in $[0,1]$. The closedness of $I$ follows easily from property (iii). Hence $I=[0,1]$.

\smallskip
\emph{Step 2: There exists a continuous map $\xi\colon B\to M$ with $\phi\compose\xi = \id_B$.} This is a standard monodromy argument: By property (ii) in the quasicovering definition, there exists a point $z_0\in M$; let $x_0=\phi(z_0)$. Every point $x_1\in B$ can be connected to $x_0$ by a path $\gamma$, and Step 1 yields a unique path $\tilde{\gamma}$ in $M$ with $\phi\compose\tilde{\gamma}=\gamma$ and $\tilde{\gamma}(0)=z_0$. We have to prove that $\xi(x_1)\define\tilde{\gamma}(1)$ does not depend on the choice of $\gamma$. This follows from the simply-connectedness of $B$, because it is straightforward to verify that homotopic choices of $\gamma$ yield the same $\tilde{\gamma}(1)$. It remains to check that the resulting map $\xi\colon B\to M$ is continuous, which is also straightforward. (Cf.\ e.g.\ \cite[proof of Theorem III.4.1]{Bredon}.)

\smallskip
\emph{Step 3: $\xi\compose\phi = \id_M$ holds.} The set $S\define\set{z\in M\suchthat \xi(\phi(z))=z}$ is nonempty because it contains $z_0$.

\smallskip
Let $z\in M$. There exists an open neighbourhood $U_0$ of $z$ in $M$ such that $\phi\restrict U_0$ is a diffeomorphism onto its image. There exists an open neighbourhood $U_1$ of $\xi(\phi(z))$ in $M$ such that $\phi\restrict U_1$ is a diffeomorphism onto its image. Since $W'\define\phi(U_0)\cap\phi(U_1)$ is a neighbourhood of $\phi(z)=\phi(\xi(\phi(z)))$ in $B$, there exists a connected open neighbourhood $W$ of $\phi(z)$ whose closure in $B$ is contained in $W'$. The sets $V_i\define(\phi\restrict U_i)^{-1}(W)$ are nonempty, connected, and open in $\phi^{-1}(W)$. They are also closed in $\phi^{-1}(W)$: the closure of $V_i$ in $M$ is contained in $(\phi\restrict U_i)^{-1}(W')$, and we have $(\phi\restrict U_i)^{-1}(W')\cap \phi^{-1}(W) = (\phi\restrict U_i)^{-1}(W)$. Thus $V_0$ and $V_1$ are connected components of the manifold $\phi^{-1}(W)$, hence either equal or disjoint.

\smallskip
The set $V\define V_0\cap(\xi\compose\phi)^{-1}(V_1)$ is an open neighbourhood of $z$ in $M$. If $x=\xi(\phi(x))$ holds for some $x\in V$, then $\xi(\phi(x))\in V_0\cap V_1$ and thus $V_0=V_1$. In that case $y=\xi(\phi(y))$ holds for every $y\in V$: the points $y$ and $\xi(\phi(y))$ lie both in $V_1$ and have the same $\phi$-image, and $\phi\restrict V_1$ is injective.

\smallskip
Therefore $S$ and $M\without S$ are open in $M$: if one of these sets contains $z$, then it contains the neighbourhood $V$ of $z$. Since $M$ is connected, we obtain $S=M$. This completes the proof of Step 3.

\medskip
The steps 2 and 3 show that $\phi$ is a homeomorphism. Since it is a local diffeomorphism, it is a diffeomorphism.
\end{proof}


\section{A proposition}

Recall that a map $f\colon M\to N$ from a manifold $M$ to a Lorentzian manifold $(N,h)$ is \emph{spacelike} iff for every $x\in M$ the image of $T_xf \colon T_xM\to T_{f(x)}N$ is spacelike; here the subspace $\set{0}$ of $T_{f(x)}N$ counts as spacelike.

\begin{lemma} \label{pathlemma}
Let $n\geq0$ and $c\in\R_{\leq0}$, let $w\colon\cointerval{0}{1}\to\M^{n,1}_c$ be a spacelike path such that $\pr\compose w\colon\cointerval{0}{1}\to\R^n$ has finite euclidean length. Then $w$ has finite length.
\end{lemma}
\begin{proof}
For $y\in\M^{n,1}_c=\R^n\times\R$, the map $T_y\pr\colon T_y\M^{n,1}_c=\R^n\times\R \to T_{\pr(y)}\R^n=\R^n$ is given by $(u,w)\mapsto u$. We claim that $\abs{v}_{\M^{n,1}_c} \leq \abs{(T_y\pr)(v)}_\eucl$ holds for all $\M^{n,1}_c$-spacelike $v$. This is obvious for $c=0$: $\abs{(u,w)}_{\M^{n,1}_0}^2 = \abs{u}_\eucl^2 -w^2 \leq \abs{u}_\eucl^2 = \abs{(T_y\pr)(u,w)}_\eucl^2$. For $c<0$, we have $\abs{(u,w)}_{\M^{n,1}_c}^2 = g_{n,2}\big(T_y\varpi(u,w),T_y\varpi(u,w)\big)$ (cf.\ Notation \ref{notation}), where $T_y\varpi(u,w)\in T_{\varpi(y)}\mathcal{H}^{n,1}_c\subseteq \R^n\times\R^2$ has the form $(u,b(y,u,w))$ for some $b(y,u,w)\in\R^2$. Thus $\abs{(u,w)}_{\M^{n,1}_c}^2 = \abs{u}_\eucl^2 -\abs{b(y,u,w)}_\eucl^2 \leq \abs{u}_\eucl^2 = \abs{(T_y\pr)(u,w)}_\eucl^2$. This proves our claim.

\smallskip
We obtain $\length(w) = \int_0^1\abs{w'(t)}\diff t \leq \int_0^1\bigabs{T_{w(t)}\pr(w'(t))}_\eucl\diff t
= \int_0^1\abs{(\pr\compose w)'(t)}_\eucl\diff t = \length(\pr\compose w)$.
\end{proof}

We say that a map $f\colon(M,g)\to(N,h)$ from a Riemannian manifold to a Lorentzian manifold is \emph{long} iff it is spacelike and for every interval $I\subseteq\R$ and every path $w\colon I\to M$, the $g$-length of $w$ is finite if the $h$-length of $f\compose w$ is finite. For example, every spacelike isometric immersion is long.

\begin{proposition} \label{rigidityprp}
Let $n\geq0$ and $c\in\R_{\leq0}$, let $(M,g)$ be a nonempty connected complete Riemannian $n$-manifold, let $f\colon(M,g)\to\M^{n,1}_c$ be a long immersion. Then $f$ is a smooth imbedding, and $\pr\compose f\colon M\to\R^n$ is a diffeomorphism.
\end{proposition}
\begin{proof}
The map $\phi\define\pr\compose f$ is an immersion, because for every $x\in M$ the image of $T_xf\colon T_xM\to T_{f(x)}\M^{n,1}_c$ is spacelike and $T_{f(x)}\pr$ maps every spacelike subspace of $T_{f(x)}\M^{n,1}_c$ injectively to $T_{\pr(f(x))}\R^n$ (since $\ker(T_{f(x)}\pr) = \set{0}\times\R \subseteq \R^n\times\R=T_{f(x)}\M^{n,1}_c$ is timelike). We claim that $\phi$ is a quasicovering.

\smallskip
Let $\gamma\colon[0,1]\to\R^n$ and $\tilde{\gamma}\colon\cointerval{0}{1}\to M$ be paths with $\phi\compose\tilde{\gamma}=\gamma\restrict\cointerval{0}{1}$. The path $\pr\compose f \compose\tilde{\gamma} = \gamma\restrict\cointerval{0}{1}$ in $\R^n$ has finite euclidean length because $\gamma$ has finite euclidean length. By Lemma \ref{pathlemma}, $f\compose\tilde{\gamma}$ has finite length. Since $f$ is long, $\tilde{\gamma}$ has finite $g$-length.

\smallskip
We choose a sequence $(t_k)_{k\in\N}$ in $\cointerval{0}{1}$ which converges to $1$. Since $\tilde{\gamma}$ has finite $g$-length, there is no $\eps>0$ such that $\forall k_0\in\N \colon \exists k,l\geq k_0 \colon \dist_g(\tilde{\gamma}(t_k),\tilde{\gamma}(t_l))\geq\eps$. Thus $(\tilde{\gamma}(t_k))_{k\in\N}$ is a Cauchy sequence in $(M,g)$. Completeness implies that it converges to some point $x\in M$. We extend $\tilde{\gamma}$ to $[0,1]$ by $\tilde{\gamma}(1)=x$. Using that $\phi$ maps a neighbourhood of $x\in M$ diffeomorphically to its image, we obtain $\phi(\tilde{\gamma}(1)) = \phi(\lim_{k\to\infty}\tilde{\gamma}(t_k)) = \lim_{k\to\infty}\phi(\tilde{\gamma}(t_k)) = \lim_{k\to\infty}\gamma(t_k) =\gamma(1)$ and deduce the smoothness of the extended $\tilde{\gamma}$ from $\gamma = \phi\compose\tilde{\gamma}$.

\smallskip
This shows that $\phi$ is a quasicovering, as claimed. By Lemma \ref{lemmaquasicovering}, $\phi$ is a diffeomorphism. Since $\phi$ is injective, so is $f$.

\smallskip
Moreover, $f$ is proper, i.e., $f^{-1}(C)$ is compact for every compact set $C\subseteq\M^{n,1}_c$. That's because $\pr(C)$ and thus $(\pr\compose f)^{-1}(\pr(C))$ are compact and $f^{-1}(C)$ is a closed subset of $(\pr\compose f)^{-1}(\pr(C))$.

\smallskip
Since every proper injective immersion is a smooth imbedding, the proof is complete.
\end{proof}

\emph{Remark.} We will apply Proposition \ref{rigidityprp} only in a situation where we know already that $M$ is simply connected. But that information would not simplify the proof.


\section{Proof of Theorem \ref{rigid}}

\begin{lemma} \label{toplemma}
Let $n\geq0$, let $M$ be a connected $n$-manifold which contains a simply connected noncompact $n$-di\-mensional sub\-mani\-fold-with-boundary that is closed in $M$ and has compact boundary. Then every covering map $\pi\colon\R^n\to M$ is a diffeomorphism.
\end{lemma}
\begin{proof}
When a connected $1$-manifold $M$ contains a noncompact subset which is closed in $M$, then $M$ is diffeomorphic to $\R$. Thus the lemma is true for $n=1$. The case $n=0$ is even simpler. Now we assume $n\geq2$. Let $Z$ be a simply connected noncompact $n$-sub\-mani\-fold-with-boundary of $M$ which is closed in $M$ and has compact boundary.

\smallskip
Since $Z$ is simply connected, the submanifold-with-boundary $\pi^{-1}(Z)$ of $\R^n$ is the disjoint union of connected components $\tilde{Z}_i$ such that $\pi\restrict\tilde{Z}_i\colon\tilde{Z}_i\to Z$ is a diffeomorphism. In particular, each $\tilde{Z}_i$ has compact boundary. Thus the boundary of $\pi^{-1}(Z)$ is a disjoint union of countably many compact nonempty connected $(n-1)$-manifolds $\Sigma_j$. No connected component $\tilde{Z}_i$ of $\pi^{-1}(Z)$ is compact, because otherwise $\pi(\tilde{Z}_i)=Z$ would be compact.

\smallskip
For each $j$, the Jordan/Brouwer separation theorem (cf.\ \cite{Lima1988} for a simple proof) implies that $\R^n\without\Sigma_j$ has precisely two connected components. Since $n\geq2$, precisely one of these two components is relatively compact in $\R^n$ (namely the unique component whose closure in the one-point compactification $S^n=\R^n\cup\set{\infty}$ of $\R^n$ does not contain the point $\infty$); we call it $\interior_j$ and denote the closure of the other component by $\exterior_j$.

\smallskip
We claim that for each $j$, \;$\pi^{-1}(Z)$ is contained in $\exterior_j$. Assume not. Then $\pi^{-1}(Z)\cap\interior_j\neq\leer$. Either a connected component of $\pi^{-1}(Z)$ is contained in $\interior_j$, or $\pi^{-1}(Z)$ touches $\Sigma_j$ from the interior (that is, $U\cap\interior_j\cap\pi^{-1}(Z)\neq\leer$ holds for every neighbourhood $U$ of $\Sigma_j$ in $\R^n$). Since $\Sigma_j$ is a boundary component of $\pi^{-1}(Z)$, the latter alternative implies that $\Sigma_j$ has a neighbourhood $U$ with $U\cap(\exterior_j\without\mfbd\,\exterior_j)\cap\pi^{-1}(Z)=\leer$. In each case, there exists a connected component $\tilde{Z}_i$ of $\pi^{-1}(Z)$ which is contained in the closure of $\interior_j$. Since $\pi^{-1}(Z)$ is closed in $\R^n$ (because $Z$ is closed in $M$), this $\tilde{Z}_i$ is compact. This contradiction proves our claim.

\smallskip
Thus $\pi^{-1}(Z)$ is contained in $\bigcap_j\exterior_j$ (which is by definition equal to $\R^n$ if the index set is empty). The two sets are even equal, for otherwise a boundary component $\Sigma_j$ of $\pi^{-1}(Z)$ would meet the interior of $\bigcap_j\exterior_j$, which is not possible because $\Sigma_j=\mfbd\,\exterior_j$ is contained in the boundary of $\bigcap_j\exterior_j$.

\smallskip
We claim that $\bigcap_j\exterior_j$ is connected. To show this, consider $x,y\in \bigcap_j\exterior_j$. We modify the straight path $\gamma$ in $\R^n$ from $x$ to $y$ on each interval $[a,b]$ it spends in $\interior_j$ for some $j$: since $\gamma(a),\gamma(b)$ lie in $\Sigma_j$, we can replace $\gamma\restrict[a,b]$ by a path in $\Sigma_j$ from $\gamma(a)$ to $\gamma(b)$. This yields a path from $x$ to $y$ in $\bigcap_j\exterior_j$ and thus proves our claim.

\smallskip
Hence $\pi^{-1}(Z)$ is connected, and $\pi$ maps $\pi^{-1}(Z)$ diffeomorphically to $Z$. The connectedness of $M$ implies that $\pi$ is a one-sheeted covering, i.e.\ a diffeomorphism.
\end{proof}

\emph{Remark.} In applications to positive energy theorems, one has much more information than is assumed in Lemma \ref{toplemma}: one knows that $M$ (of dimension $n\geq3$) is noncompact and contains a compact $n$-dimensional submanifold-with-boundary $C$ such that each connected component $Y$ of $M\without C$ is diffeomorphic to $S^{n-1}\times\oointerval{0}{1}$; the closure $Z$ in $M$ of each of these ends $Y$ is a submanifold-with-boundary of $M$ which is diffeomorphic to $S^{n-1}\times\cointerval{0}{1}$ and thus satisfies the assumptions of the lemma. But all this additional information would not help much in the proof. For instance, $\pi^{-1}(C)$ could a priori still be noncompact; this makes arguments involving ends difficult.

\begin{proof}[Proof of Theorem \ref{rigid}]
Let $\pi\colon\tilde{M}\to M$ be the universal covering of $M$, let $\tilde{g}\define \pi^\ast g$, let $\tilde{N}$ be the pullback bundle $\pi^\ast N$ over $\tilde{M}$, and define $\tilde{K}=\pi^\ast K\in\Gamma(\Sym^2T^\ast\tilde{M}\otimes\tilde{N})$ by $\tilde{K}(v,w)=K(\pi_\ast v,\pi_\ast w)\in N_{\pi(x)}=(\pi^\ast N)_x$ for all $x\in\tilde{M}$ and $v,w\in T_x\tilde{M}$. Since $(M,g,N,K)$ satisfies the Gauss and Codazzi equations for constant curvature $c$, so does $(\tilde{M},\tilde{g},\tilde{N},\tilde{K})$. Being the pullback of a complete metric by a covering map, $\tilde{g}$ is complete.

\smallskip
Proposition \ref{bgm} tells us that there exists an isometric immersion $(f,\iota)$ of $(\tilde{M},\tilde{g},\tilde{N},\tilde{K})$ into $\M^{n,1}_c$; and that any two such immersions differ by an isometry of $\M^{n,1}_c$. Proposition \ref{rigidityprp} implies that $f$ is an isometric imbedding and that $\pr\compose f\colon\tilde{M}\to\R^n$ is a diffeomorphism. We identify $\tilde{M}$ with $\R^n$ via $\pr\compose f$.

\smallskip
Lemma \ref{toplemma} shows that the covering $\pi\colon\R^n\to M$ is a diffeomorphism. $(\tilde{M},\tilde{g},\tilde{N},\tilde{K})$ and $(M,g,N,K)$ can be identified via $\pi$, and the theorem follows.
\end{proof}

\medskip
\emph{Remark 1.} The proof here is similar to the work of Maerten \cite[second half of the proof of the first theorem in Section 4]{Maerten2006} (which deals with the case $c<0$ on a spin manifold) insofar as both employ the universal covering of $M$ and argue that it is one-sheeted. Maerten uses apparently a statement similar to Lemma \ref{toplemma} at the end of his proof, but does not give a reference or spell out the details.

\medskip
\emph{Remark 2.} The proof of the positive energy theorem in \cite{ParkerTaubes} yields already the information that the hypersurface $M$ has only one end in the rigidity case. The arguments above provide a second, independent proof that $M$ has only one end.


\bigskip
\thanks{\emph{Acknowledgement.} I would like to thank Olaf Müller for a helpful discussion.}



\newcommand{\dummysort}[1]{}\def\cprime{$'$}
  \def\Dbar{\leavevmode\lower.6ex\hbox to 0pt{\hskip-.23ex \accent"16\hss}D}
  \def\Dbar{\leavevmode\lower.6ex\hbox to 0pt{\hskip-.23ex \accent"16\hss}D}
  \def\Dbar{\leavevmode\lower.6ex\hbox to 0pt{\hskip-.23ex \accent"16\hss}D}

\end{document}